\documentclass{article}

\usepackage[utf8]{inputenc}   
\usepackage[T1]{fontenc}      
\usepackage[english]{babel}   

\usepackage{amsmath, amsfonts, amssymb}  

\usepackage{graphicx}         

\usepackage{ragged2e}         

\usepackage[hidelinks]{hyperref}  

\title{THE EXTENDED ALPHA GROUP DYNAMIC MAPPING}

\author{
Cleber Souza Corr\^{e}a\\
Instituto de Aeron\'{a}utica e Espa\c{c}o, S\~{a}o Jos\'{e} dos Campos, SP, Brazil\\
\href{mailto:clebercsc@fab.mil.br}{E-mail: clebercsc@fab.mil.br}\\ 
\and
Thiago Braido Nogueira de Melo\\
Instituto de Aeron\'{a}utica e Espa\c{c}o, S\~{a}o Jos\'{e} dos Campos, SP, Brazil\\
\href{mailto:braidotbnm@fab.mil.br}{E-mail: braidotbnm@fab.mil.br}\\
}

\begin{document}

\maketitle
\selectlanguage{english}

\justifying

\noindent \textbf{Mathematics Subject Classification (2020):} 34Cxx, 34Dxx, 37Cxx, 37Nxx, 37Mxx.

\begin{abstract}
This paper investigates the qualitative behavior of a system of ordinary differential equations (ODEs) defined by a matrix operator derived from the algebraic structure of the Alpha Group. The system depends on a rotational parameter that continuously deforms the underlying geometry of the phase space.

Using a fourth-order Runge--Kutta numerical scheme, we analyze the evolution of trajectories and identify the presence of critical parameter values at which the system undergoes qualitative transitions. In particular, we observe the emergence of critical dynamical regions associated with changes in the interaction between dynamically defined subspaces. As the rotation parameter varies from $0$ to $\pi/2$, the system transitions from a regime with Euclidean-type geometric behavior to a non-Euclidean configuration induced by the Alpha Group structure.

These transitions correspond to changes in stability and global phase space organization, including the formation of invariant structures and attractor-like behavior at infinity. The results suggest that the underlying matrix operator acts as a generator of structured transformations governing the system dynamics.

This work provides a computational and qualitative framework for studying parameter-dependent dynamical systems with evolving geometric structure.
\end{abstract}

\noindent\textbf{Keywords}: Ordinary differential equations, Asymptotic Compactification, Hyperboloid, Manifolds, Hypercomplex Topology.

\section{Introduction}

The use of ordinary differential equation (ODE) systems is a common method in science to study the dynamics of complex systems. This methodology is well-established in the literature, where ODE simulations are used to observe resulting behavior and map properties by changing initial conditions. The Alpha Group defines a structure based on a division-based operation. Corr\^ea et al. (2022) proposed the Alpha group, using group theory, formed by the transformation of two infinite planes interacting through a change in the division operation between them. This interaction creates a third element with morphism and preserves the operations in both planes. Infinity is associated with a geometric representation that induces topological deformations and generates an attractor-like structure in $R^4$. This situation resulted from a division operation associated with the rotation between the planes that form the Alpha group. Corr\^ea et al. (2024) demonstrated that the metrics of infinitesimal distance from Riemannian and Euclidean space are special cases of the Alpha group's metrics. The Alpha group satisfies the properties defined by group theory. These findings establish a consistent structure that geometrically and topologically characterizes a hypercomplex numerical group, possessing inherent properties that define transformations between surfaces. These definitions enable the development of new mathematical research.

The geometry of the alpha group is based on George Cantor's theories, which explored the nature of different types of infinity (Corrêa and de Melo, 2025). This allows for the interpretation of the existence of numerous and varied types of infinity and their connection to geometry and topology. This structure has a trigonometric representation associated with the hypercomplex plane, resulting in a $4 \times 4$ antisymmetric matrix with a system that incorporates tangent and cotangent functions between its elements. The matrix's properties must be analyzed dynamically. In such a scenario, ODE systems enable changes in initialization conditions and observation of the resulting dynamic behavior. A significant feature of the Alpha Group is that when zero radian rotation occurs in the formation of the numeric group, it is associated with Euclidean topology. However, when the rotation is $\pi/2$ radians, maximum deformation occurs, generating a tangent plane at infinity. In this aspect, it converges asymptotically and may be a manifestation of asymptotic symmetry, where the field tends to a stable configuration at infinity. This stable structure is associated with a type of global attractor in four dimensions and is defined by an imaginary number ($\mu$) that has geometric and topological representations in the context of hyperbolic topology. In a geometric context, the global attractor can represent an asymptotic symmetry, where the field tends to a stable or homogeneous configuration as it approaches infinity. This topological change along the parameter $\theta$, from the local Euclidean space to the infinite tangent plane, reflects a fundamental concept that can be associated with gauge-like theory. The idea is that local symmetries can extend to describe global properties, converging asymptotically into a stable asymptotic convergence (closed hyperbolic attractor). Thus, this analysis may provide a qualitative framework for relating geometric deformations to gauge-like transformation structures. Accurately describing both local interactions and global transformations. Because the resulting matrix model exhibits properties that can change the topology, transitioning from a Euclidean topology to an $R^4$ topology in the Alpha group space, this work aims to map the dynamic characteristics of the matrix resulting from group theory, which leads to the Alpha group, using different initializations.

\section{Methodology}

\subsection{Mathematical Model}

The Extended Alpha Group can be constructed through a tensor-based division operation between two matrices formulated via the De Moivre identity, even if the quaternion in the Alpha Group can be made of two complex planes. The construction of the division approach is analogous to the Kronecker product, as demonstrated in Graham's (2018) equation (I). The result of equation (II) emerged naturally from the relationship between the operators; however, there may be other possible ways of representing this operation. Looking at the trigonometric relationship, the ratio between $\sin\theta$ and $\cos\theta$ represents the tangent ($\tan\theta$), and the cotangent ($\cot\theta$) is the inverse function of the tangent. The tangent of an angle can also be defined by the ratio between the measure of the opposite side and the measure of the adjacent side to the angle $\theta$ in a right triangle. It can replace trigonometric mathematical relations and possibly division operations. The matrix $A$ (II) resulting from the division of two complex planes generates a transformation into a matrix of a hypercomplex space. Matrix $A$ has some interesting properties:
\begin{itemize}
    \item Antisymmetry: The matrix is not symmetric, meaning that $A \neq A^T$, where $A^T$ is the transpose of $A$. This implies distinct properties compared to symmetric matrices.
    \item Non-zero determinant: The determinant of $A$ is non-zero, meaning the matrix is non-singular. This implies that the matrix has an inverse.
    \item Non-zero diagonal elements: The diagonal elements of $A$ are nonzero; the matrix has non-zero diagonal elements, but is not diagonal.
    \item Elements off the main diagonal: Elements off the main diagonal have a specific relationship based on angle $\theta$, indicating a non-trivial structure.
    \item Parameter dependency: Matrix $A$ depends on the parameter $\theta$, which implies that its properties can vary based on this parameter.
    \item Transformation properties: Matrix $A$ can represent a linear transformation in a vector space, with specific properties related to the application of that transformation.
\end{itemize}
Accordingly, the antisymmetric matrix $A$  exhibits dynamic properties as the rotation between the planes occurs to form the Alpha Group (Corr\^ea et al., 2022 and 2024). One way to analyze this matrix $A$ is to associate it with a system of ordinary differential equations (ODEs). An ordinary differential equation can be written as:
\begin{equation*}
    \frac{d}{dt}x = A \cdot x \quad (\text{I})
\end{equation*}
An improvement to matrix (II) can be achieved by inserting the parameter $\mu$ along the main diagonal. This modification characterizes the operation defined within the domain of the Alpha Group, taking into account both the effects of topological deformation between planes due to rotations in the interval $[0, 2\pi]$ and the numerical construction principles underlying the Alpha Group.

The inclusion of $\mu$ on the main diagonal ensures that the matrix is non-singular (i.e., its determinant is non-zero), while also preserving its antisymmetric structure. This property is essential for maintaining the coherence of the internal transformations governed by the group.

The matrix $M(\theta)$ is then defined as the product:
\[
M(\theta) = A(\theta) \cdot B(\mu),
\]
where $A(\theta)$ is an angular matrix explicitly dependent on the rotational parameter $\theta$, and $B(\mu)$ is a phase amplitude matrix incorporating the deformation parameter $\mu$. The resulting matrix $M(\theta)$ is a non-Hermitian structure that acts as a generator of internal vectorial variations.

In the context of the Alpha Group algebra, this composition formalizes a fundamental mechanism of internal dynamics, coupling angular and amplitude contributions into a unified framework.

$$
M(\theta) = \begin{pmatrix}
1 & -\cot\theta & -\tan\theta & 1 \\
\cot\theta & 1 & -1 & -\tan\theta \\
\tan\theta & -1 & 1 & -\cot\theta \\
1 & \tan\theta & \cot\theta & 1
\end{pmatrix}
\cdot
\begin{pmatrix}
1 & 1 & 1 & 1 \\
1 & i & 1 & 1 \\
1 & 1 & \mu & 1 \\
1 & 1 & 1 & i\mu
\end{pmatrix} \quad (\text{II})
$$

The system of ODEs is represented as:
\begin{equation*}
    \frac{d}{dt}\begin{pmatrix}
        x_1 \\
        x_2 \\
        x_3 \\
        x_4
    \end{pmatrix} = \begin{pmatrix}
        1 & -\cot\theta & -\tan\theta & 1 \\
        \cot\theta & 1 & -1 & -\tan\theta \\
        \tan\theta & -1 & 1 & -\cot\theta \\
        1 & \tan\theta & \cot\theta & 1
    \end{pmatrix}\begin{pmatrix}
        x_1 \\
        x_2 \\
        x_3 \\
        x_4
    \end{pmatrix} \quad (\text{III})
\end{equation*}
The ODE system (III) can express each row as a differential equation. Let's use variables $x_1$, $x_2$, $x_3$, and $x_4$ for the dependent variables and assume they are functions of some independent variable, usually denoted as $t$ (time). Within the ODE system with matrix $A$, $dx/dt$ is initialized as a vector $x_0 = (1, 1, 1, 1)$ with the same size as $x$. The four differential equations are then defined based on the specific model of the system: The derivative of $x_1$ for time is defined by the equation $x_1 - \cot\theta \cdot x_2 - \tan\theta \cdot x_3 + x_4$, and the other three derivatives ($x_2$, $x_3$, $x_4$) are defined similarly. Essentially, this is a system of first-order ODEs with four differential equations, each representing the rate of change of a system variable over time. This system is then solved numerically using the fourth-order Runge-Kutta method. A script was built to numerically simulate the ODE system, written in Python, and the simulation environment was the website mycompiler.io (https://www.mycompiler.io/pt). In this script, the fourth-order Runge-Kutta function was used. The discretization factor $h$ was $0.001$ and integrated over time at $1.5$. The value of $\mu$ was also replaced by $1$. The Python script was programmed to generate the Poincar\'e map associated with each rotation angle of the matrix $A$ system. The NumPy Python package is the fundamental package for scientific computing in Python. It is a Python library that provides a multidimensional array object, various derived objects (such as masked arrays and matrices), and an assortment of routines for fast operations on arrays, including mathematical, logical, shape manipulation, sorting, selecting, discrete Fourier transforms, basic linear algebra, basic statistical operations, random simulation, and much more.

\subsection{The Poincar\'e Map}

The Poincaré map is a fundamental technique in the study of dynamical systems, providing a way to analyze the behavior of differential equations by projecting trajectories onto a lower-dimensional plane. This method simplifies the analysis by reducing the continuous dynamics to a discrete system, allowing the identification of periodic orbits, bifurcations, and attractors.

The central idea behind the Poincaré map is to define a specific subset of the system's phase space, known as the \textit{Poincaré section}, and to observe the intersection points of the trajectories with this section. Typically, this is done by fixing one of the system's variables and recording the points at which the trajectory crosses the section transversally.

These intersection points generate a two-dimensional map—called the \textit{Poincaré map}—which captures the essence of the system's dynamics in a discrete form. It is particularly useful in the analysis of complex or chaotic systems where continuous observation is difficult or uninformative.

In computational implementations, such as in Python, a custom function (e.g., \texttt{poincare\_map(x\_values)}) can be defined to generate the Poincaré map numerically. This function typically extracts the coordinates of the points where the trajectories intersect the chosen section, using tools from numerical integration and array handling (e.g., NumPy).

This approach allows for an effective visualization and classification of dynamic behavior, serving as a bridge between continuous differential systems and discrete dynamics.

\subsection{Phase Diagram}

The Python script also generated the phase diagram, which graphically represents the system’s trajectories within the state space defined by a system of four differential equations. These trajectories illustrate how the state variables evolve and interact dynamically.

In a phase diagram, the evolution of time is visualized through projections in selected pairs of state variables, allowing a two-dimensional analysis of the system’s behavior. This method provides valuable insight into the structure and stability of the system.

The interpretation of the phase diagram depends on the specific characteristics of the dynamic system to be studied. Among the general features that may appear in such diagrams are:
\begin{itemize}
    \item Trajectories and their global behavior over time;
    \item Equilibrium points (fixed points);
    \item Attractors and repellors;
    \item Stability regions and their boundaries;
    \item Limit cycles and periodic orbits;
    \item Saddle points and separatrices;
    \item Emergent chaotic behavior.
\end{itemize}

Phase diagrams are essential tools for identifying and classifying dynamic regimes, especially in systems where analytical solutions are difficult or impossible. When combined with Lyapunov analysis and Poincaré maps, they offer a powerful framework for understanding both local and global stability properties.

\subsection{Lyapunov Function}

To identify attractors at infinity, it is often necessary to employ specific numerical techniques capable of analyzing the long-term behavior of a dynamical system. In this context, the Lyapunov method was used to assess both the stability and the presence of attractors located at infinity.

The key idea is that if a system possesses an attractor at infinity, the state variables tend to grow without bound, often at an exponential rate. This growth can be quantified using a Lyapunov function evaluated along the solution trajectories of the system.

Numerically, the Lyapunov values were computed using the following function:

\begin{verbatim}
lyapunov_values = np.apply_along_axis(lyapunov_function, 1, ode_solution)
\end{verbatim}

where \texttt{ode\_solution} is obtained by numerical integration of the system using:

\begin{verbatim}
from scipy.integrate import odeint
\end{verbatim}

The \texttt{SciPy} library is a collection of mathematical algorithms and convenience functions built on top of \texttt{NumPy}. It adds substantial computational power to Python by offering high-level routines and data manipulation tools, making it suitable for tasks such as the numerical analysis of dynamical systems, including the detection of asymptotic behaviors and attractors.

This approach enables the classification of system stability through the sign and magnitude of the Lyapunov values, especially when investigating structures that diverge toward infinity or stabilize on non-trivial asymptotic manifolds.

\subsection{Bifurcation Diagram}

The bifurcation diagram of the differential equation system defined by matrix $A$ was generated using the Python libraries \texttt{NumPy} and \texttt{matplotlib.pyplot}. The \texttt{NumPy} library plays a fundamental role in scientific computing with Python, providing efficient support for multidimensional arrays and high-performance mathematical operations.

To solve the system, a custom implementation of the classical fourth-order Runge–Kutta method was employed. This method numerically integrates the system's dynamics over a range of initial conditions or parameters, focusing on the variation concerning the angular parameter $\theta$.

The function calls:

\begin{verbatim}
bifurcation_data.append(solution_rk4[-1])
\end{verbatim}

was used to collect the final value of the state variables from each Runge–Kutta integration. These terminal values were then plotted to construct the bifurcation diagram, which reveals changes in the qualitative behavior of the system as $\theta$ varies.

The resulting diagram provides insight into the structure of the system's dynamics, including the appearance of fixed points, periodic solutions, branching behavior, and transitions that may indicate bifurcations or the onset of chaotic regimes.

\subsection{Jacobian Matrix}

The change in the ODE system of the matrix $A$ is associated with the $\theta$ angle. If the maximum deformation occurs at $\pi/2$ radians and the formation of an Alpha group number space results in a significant change in the system's behavior. This may indicate a critical point or singularity in the system's configuration space. The eigenvalues of the Jacobian matrix of the ODE system of matrix $A$ were also calculated close to $\pi/2$ to observe whether the majority of the eigenvalues were composed of complex numbers with a negative real part. This suggests that the system has asymptotically stable behavior. The Jacobian matrix of a system of differential equations is a matrix where each element $J_{ij}$ represents the partial derivative of the function $f_i$ to the variable $x_j$. The formula for the Jacobian matrix of the ODE system is given by:
\begin{equation*}
    J_{4x4} = \begin{pmatrix}
        \frac{\partial f_1}{\partial x_1} & \frac{\partial f_1}{\partial x_2} & \frac{\partial f_1}{\partial x_3} & \frac{\partial f_1}{\partial x_4} \\
        \frac{\partial f_2}{\partial x_1} & \frac{\partial f_2}{\partial x_2} & \frac{\partial f_2}{\partial x_3} & \frac{\partial f_2}{\partial x_4} \\
        \frac{\partial f_3}{\partial x_1} & \frac{\partial f_3}{\partial x_2} & \frac{\partial f_3}{\partial x_3} & \frac{\partial f_3}{\partial x_4} \\
        \frac{\partial f_4}{\partial x_1} & \frac{\partial f_4}{\partial x_2} & \frac{\partial f_4}{\partial x_3} & \frac{\partial f_4}{\partial x_4}
    \end{pmatrix} \quad (\text{IV})
\end{equation*}
where $f_i$ is the $i$-th function of the ODE system of matrix $A$ and $x_j$ is the $j$-th state variable. The `sympy` library was used to calculate the eigenvalues of the Jacobian matrix.

\section{Results}

The results were important because they show the dynamic behavior of the system of ordinary differential equations associated with the matrix $A$. Figure 1 illustrates the system dynamics under zero-radian rotation, representing the Euclidean regime. In this aspect, its topology can be associated with an Euclidean space. The system converges to a stable zero-equilibrium state when the rotation approaches $0 + n\pi$ radians. In a certain aspect, it generates a type of dynamically stable and convergent behavior.

\begin{figure}[ht!]
    \centering
    \includegraphics[width=1.0\textwidth]{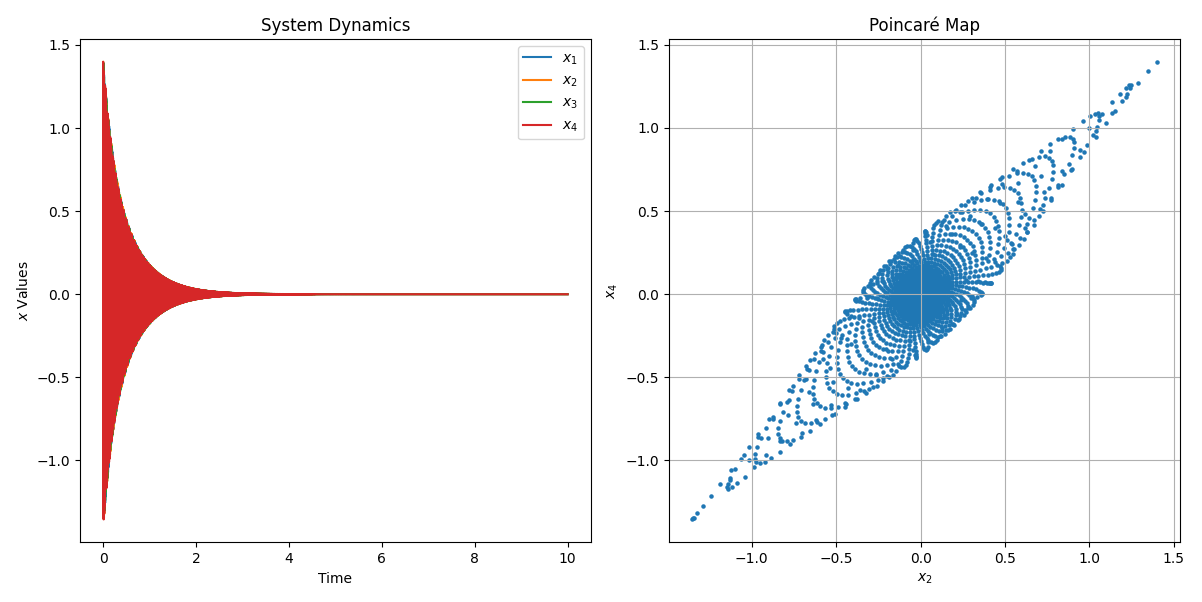}
    \caption{The dynamics of the ODE system and Poincar\'e map associated with simulation near zero radian rotation or $0 + n\pi$ radians.}
    \label{fig:fig1}
\end{figure}

An analysis of the Poincar\'e map can provide valuable information about the dynamics of a dynamical system. Figure 1 shows a distribution whose dynamic system has a stable equilibrium point. It can be seen that the points in the Poincar\'e figure converge to this point. The trajectories in the phase diagram represent the solutions of the system over time. Each point on the trajectory corresponds to the state of the system at a given moment. As shown in Figure 2, trajectory convergence behavior is shown for a specific region at zero value. The topology of Euclidean space is characterized by zero rotation between planes. It was also calculated for the angle of 0.2 degrees (Euclidean Space); however, for the zero-degree angle, the ODE system exhibits trajectories converging toward the origin, and the Lyapunov function does not grow exponentially over time but instead approaches a constant value.

\begin{figure}[ht!]
    \centering
    \includegraphics[width=1.0\textwidth]{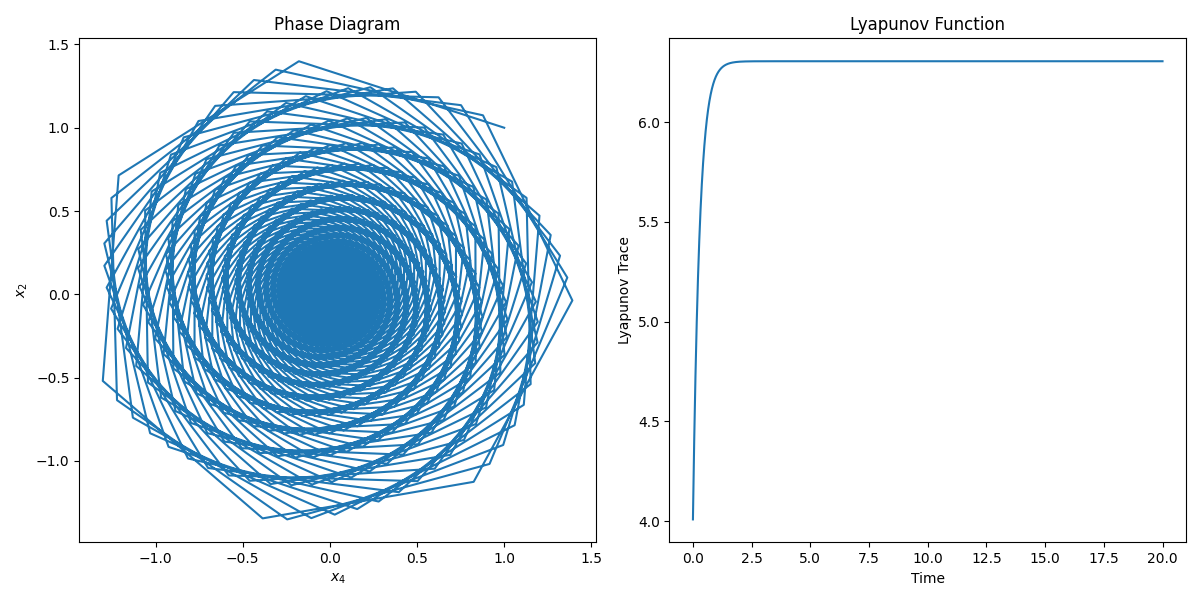}
    \caption{The phase diagram and the Lyapunov function for the Matrix A ODE system for the angle 0.02 degrees close to 0 radians.}
    \label{fig:fig2}
\end{figure}

Figure 3 shows the situation when the rotation is $\pi/2 + n\pi$ radians (90 degrees). In this case, the ODE system shows that the system grows infinitely, tending to infinity asymptotically. In this case, the maximum topological deformation and the formation of the Alpha group space are associated. In this case, the ODE with matrix $A$ can be a system behavior generator in $R^4$, within a gauge-like interpretative framework, there is a relationship between local behavior (zero radians) and global behavior (like the asymptotic attractor ($\pi/2$ radians)), which is similar to the idea of its dynamical system growing and forming complex geometric structures. The alpha group and its matrix $M(\theta)$ could be related to structured transformation mechanisms under local or global symmetry-like operations. The matrix may be interpreted as inducing structured symmetry-like and asymmetry-like transformations. The deformations are associated with the canonical vector imaginary number $\mu$ (Corr\^ea et al., 2022 and 2024). The vector element $\mu$, considered as a canonical vector and topological invariant, transcends the traditional role of an imaginary unit. Within the algebraic formulation of the Alpha Group, it $\mu$ acquires an internal vectorial condition, functioning as a generator of symmetry that organizes and structures the internal algebraic relations.

Importantly, the element $\mu$ remains unaffected by fluctuations in the angular parameter $\theta$, internal rotations governed by $M(\theta)$, or by any external perturbations. It thus becomes an absolute topological reference, serving as an invariant anchor within a complex and topologically structured vector space.

This framework may be interpreted as a transition of topological domains from a collapsed 3-sphere ($S^3$), which corresponds to a real regime, to a 4-sphere-like regime ($S^4$), representing a complex regime. In this transition, a hidden degree of freedom begins to emerge through the imaginary component of the structure. The dynamic behavior of the complex eigenvalues of $M(\theta)$ plays a central role: when $\operatorname{Re}(\lambda) \rightarrow 0$ and $\operatorname{Im}(\lambda) \rightarrow \infty$, phenomena such as may indicate torsion-like or curvature-like effects.

In this context, the vector $\mu$ acts as an idempotent and invariant element, ensuring internal coherence and stability across topological deformations. Its presence is essential for preserving the structural integrity and consistency of the Alpha Group algebra.
In this situation, it generates a type of attractor that is dynamically stable but converges asymptotically at infinity. 

\begin{figure}[ht!]
    \centering
    \includegraphics[width=1.0\textwidth]{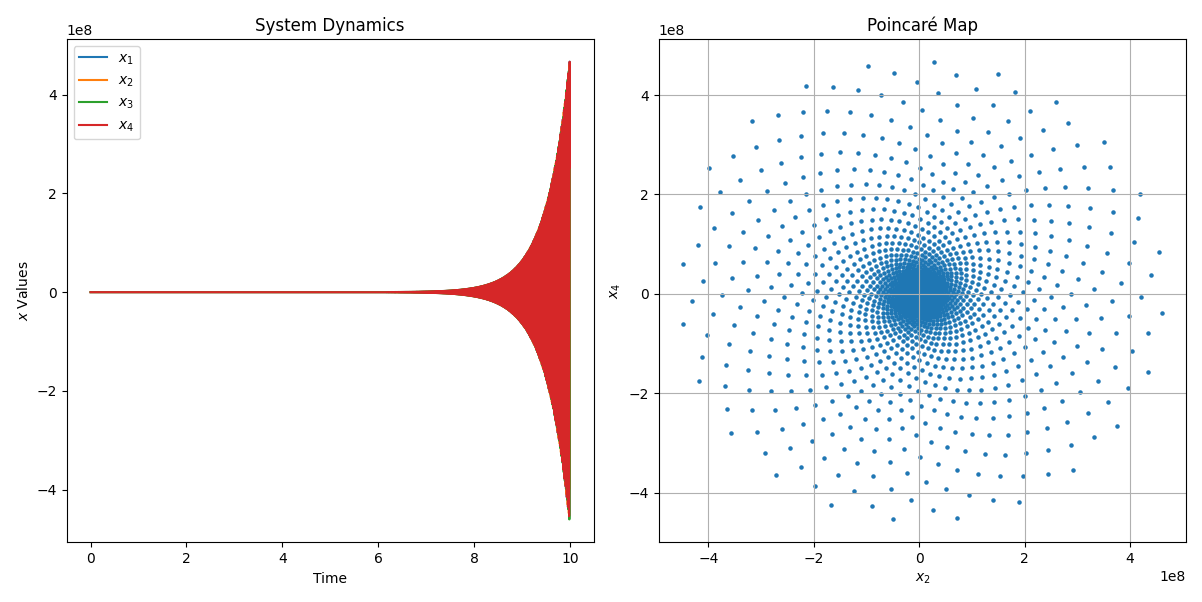}
    \caption{The dynamics of the ODE system and Poincar\'e map associated with simulation close to $\pi/2$ radians or $\pi/2 + n\pi$ radians in 10 time.}
    \label{fig:fig3}
\end{figure}

\begin{figure}[h!]
    \centering
    \includegraphics[width=1.0\textwidth]{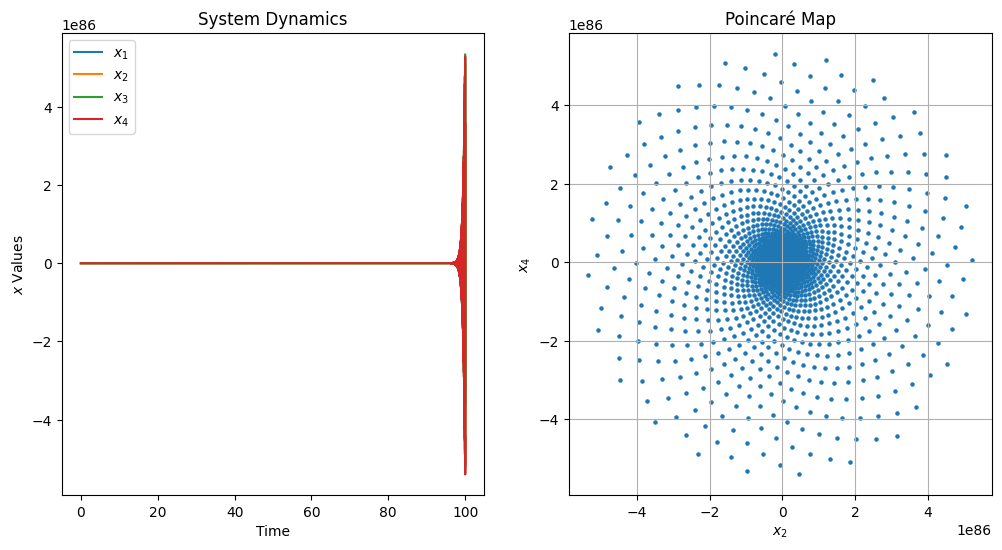}
    \caption{The dynamics of the ODE system and Poincar\'e map associated with simulation close to $\pi/2$ radians or $\pi/2 + n\pi$ radians in 100 times.}
    \label{fig:fig4}
\end{figure}

Figures 3 and 4 show that the dynamic system can be characterized as having an attractor. It can be observed in regions in the Poincaré figure where the points are concentrated, indicating that the system spends more time in these regions. Different integrations were carried out over time, and at a time of the order of $t = 100$ (Figure 4), the dynamics of the ODE system tend to present very large values and form a similar figure to the time of $t = 10$, but with orbits tending to infinity. Such dynamic behavior in $\pi/2$ indicates that in the dynamic process of the ODE system going to asymptotically infinity, the area associated with the attractor image on the Poincaré map grows proportionally to the growth rate of the ODE system values of the simulated matrix $A$. The Alpha group presents the trajectory convergence behavior for a specific region, indicating the presence of an attractor at infinity. If the phase diagram shows a convergence towards asymptotic infinity, this may be indicative of an interesting behavior of the dynamical system. The concept of a "tangent space in infinity" in closed hyperbolic topology also refers to certain notions of Gauge's theory, which often deals with how fields behave in infinity or asymptotic limits. In Gauge-like theory, there is often a relationship between local behavior (such as rotations in your system) and global behavior (as the asymptotic attractor), which is similar to the idea of its dynamic system growing and forming complex geometric structures. Therefore, this development of a tangent space at infinity and the symmetries associated with the attractor in closed hyperbolic space can be interpreted under the framework of gauge-like theory. It offers tools to understand how local transformations, such as rotations, impact the geometry and global dynamics of the system, allowing the formation of structures such as the attractor in the Alpha group. According to research by Corr\^ea et al. (2022), the Alpha group is formed by the transformation of two infinite planes that interact by changing the division operation between them, resulting in the creation of a third element with morphism and the preservation of operations in both. Infinity has a geometric representation that causes topological deformations. The geometry and topology produced by this rotation thus demonstrate the presence of a particular tensor metric in this alpha group, which has an antisymmetric and mirrored topology. Matrix $A$ has dynamic properties that can be observed in simulations carried out with the different rotations associated with initialization. Figure 5 shows the topology of the Alpha group when the rotation angle is $\pi/2 + n\pi$ radians. The observed dynamic behavior presents qualitative similarities with gauge-like transformation frameworks, especially in the context of dynamic systems and differential equations that display rotational symmetries and asymptotic behavior. Gauge-like frameworks are often associated with systems exhibiting invariance under structured transformations (such as rotations or translations), and the behavior of their system under a rotation of $\pi/2$ suggests the presence of a symmetry that is fundamental to its dynamics. Creating a type of attractor in $R^4$ is associated with the tangent plane in $\pi/2 + n\pi$ radians. The Lyapunov function (Figure 5) grows exponentially with time; this may indicate the presence of attractors at infinity. The matrix $A$ has $\mu$ on the main diagonal and possibly behaves similarly to an imaginary number. This may indicate that certain aspects of the system dynamics are strongly influenced by complex or rotating properties in the 4-dimensional state space.

\begin{figure}[ht!]
    \centering
    \includegraphics[width=1.0\textwidth]{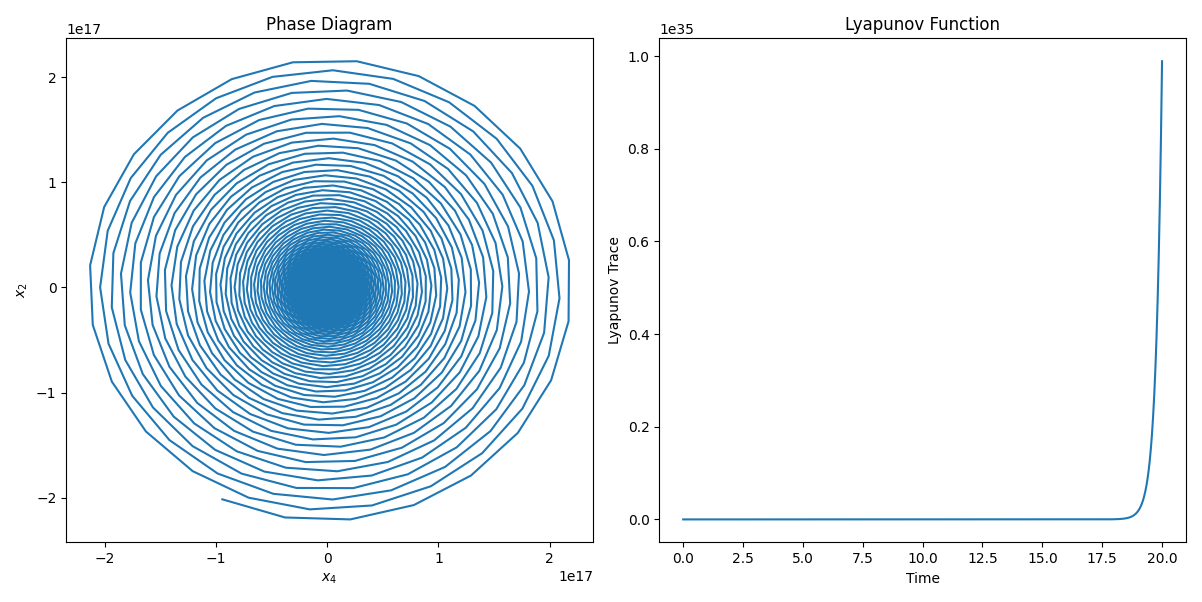}
    \caption{The phase diagram and the Lyapunov Function for the Matrix A ODE system for the angle 89.8 degrees close to $\pi/2$ radians.}
    \label{fig:fig5}
\end{figure}

However, when outside these two nodes, zero and $\pi/2$, the matrix $A$ generates a stable dynamic structure, as shown in Figure 7. This dynamic behavior presented by the numerical simulation shows that there are dynamic nodes in the ODE system analyzed in matrix $A$. The Poincar\'e map has only an orbit path, which may indicate a stable periodic orbit in the dynamical system. In dynamical systems, a periodic orbit is a closed trajectory that the system follows repeatedly over time. If this orbit is stable, it will manifest itself as an orbit path on the Poincar\'e map.

\begin{figure}[ht!]
    \centering
    \includegraphics[width=1.0\textwidth]{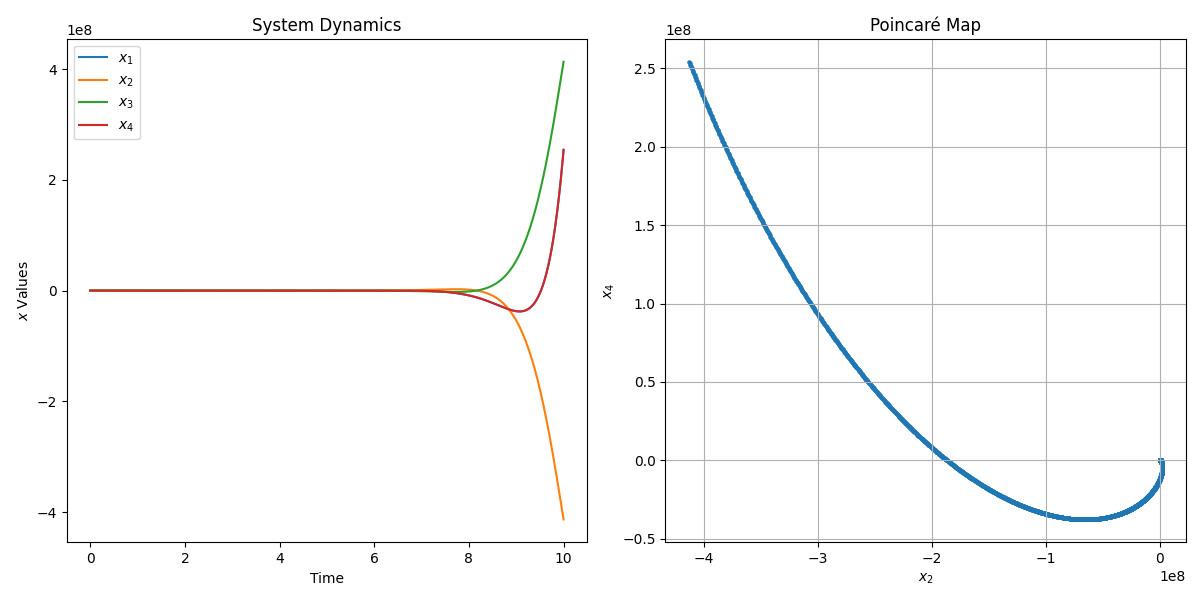}
    \caption{The dynamics of the ODE system and Poincar\'e map associated with simulation at the rotation of $\pi/6$ radians (30 degrees).}
    \label{fig:fig6}
\end{figure}

Figure 6: The phase diagram has only an orbit path, and the Lyapunov function shows the exponential curve for an angle of 30 degrees. The Lyapunov function, over time, also grows exponentially. All tested angles greater than 1 degree resulted in the Lyapunov function growing exponentially, showing a tendency to an attractor at 90 degrees.

\begin{figure}[ht!]
    \centering
    \includegraphics[width=1.0\textwidth]{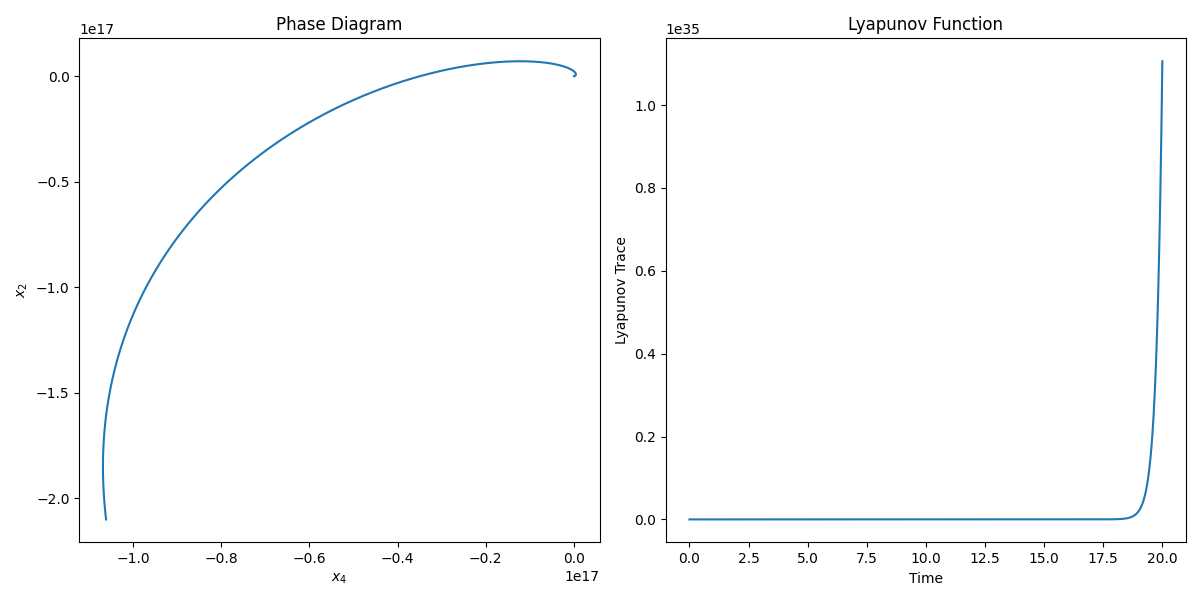}
    \caption{The phase diagram and the Lyapunov Function for the Matrix A ODE system for the angle 30 degrees close to $\pi/6$ radians.}
    \label{fig:fig7}
\end{figure}

In Figure 8, oscillatory patterns can be seen in the center of the bifurcation diagram and dispersion on the sides; this may indicate a transition to an asymptotic transformation to infinity. It was analyzed near the point $\pi/2$ using the Fast Fourier Transform (FFT), which is an efficient algorithm to calculate the Discrete Fourier Transform (DFT). DFT is an essential tool in signal processing and spectral analysis, allowing the representation of a signal in the frequency domain. The change in frequencies observed in the dynamic system indicates that there may be an interesting transition in the system's behavior when the parameter approaches $\pi/2$. The change from 40 Hz to 300–400 Hz suggests a transition to higher frequencies and may indicate a change in the dynamic regime of the system. This transition can be associated with different system behaviors. Frequency spectrum analysis is a useful tool for identifying changes in the oscillatory characteristics of a system. The presence of specific frequencies, such as 40 Hz and 300–400 Hz, may suggest the existence of periodic orbits or stable cycles in different regimes. Understanding these frequencies can provide insights into the predominant oscillation modes of this system.

\begin{figure}[ht!]
    \centering
    \includegraphics[width=1.0\textwidth]{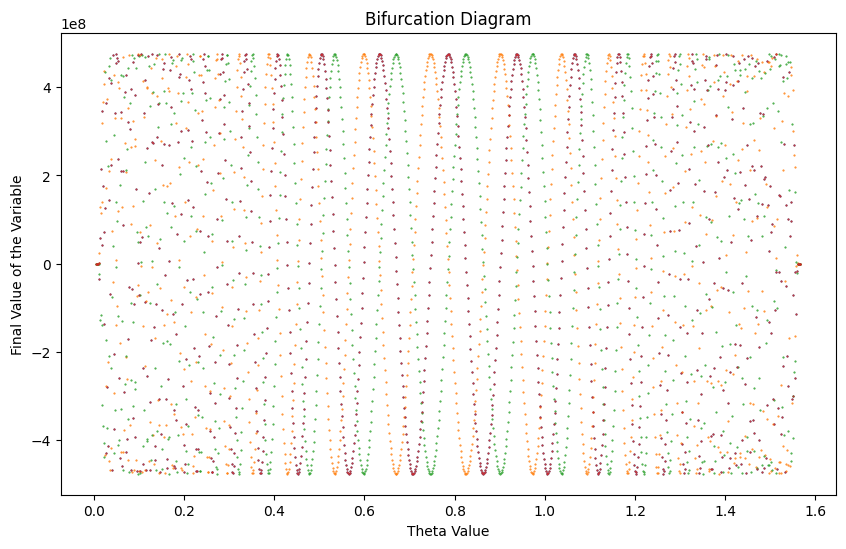}
    \caption{Bifurcation diagram associated with the angle $\theta$ and presents the dynamic behavior of the ODE system of matrix A.}
    \label{fig:fig8}
\end{figure}

The results of the eigenvalues of the Jacobian matrix in the ODE system of matrix $A$ close to $\pi/2$ are composed of complex numbers with a negative real part; this suggests that the system has asymptotically stable behavior. The result of calculating the eigenvalues of the Jacobian matrix of the ODE system of matrix $A$ was:
\begin{align*}
    \lambda_1 &= -\frac{5}{4} - \frac{1}{4} - i\frac{5}{8} - \frac{5}{8} \\
    \lambda_2 &= -\frac{5}{4} - \frac{1}{4} + i\frac{5}{8} - \frac{5}{8} \\
    \lambda_3 &= -\frac{1}{4} + \frac{5}{4} - i\frac{5}{8} + \frac{5}{8} \\
    \lambda_4 &= -\frac{1}{4} + \frac{5}{4} + i\frac{5}{8} + \frac{5}{8}
\end{align*}
In this context, stability is a desirable characteristic to ensure that the system's solutions converge to an equilibrium point. If half of the real eigenvalues are negative, this suggests that the dynamical system has stable dynamics, at least around the specific point analyzed ($\pi/2$ radians). In dynamic systems, the eigenvalues of the Jacobian matrix determine the stability of the equilibrium point. The fact that it converges towards infinity means that the system's solutions approach a stable state over time. Asymptotic stability is a desirable property, as it means that perturbations in the initial conditions or the system will be dampened over time, and the system will reach a stable state. The fact that the real eigenvalues have a negative real part suggests that any disturbance in the initial conditions will lead the system to converge to the equilibrium point at infinity. The Alpha Group existence (imaginary number $\mu$) leads to a global attractor-like state in a closed hyperbolic-type geometry, where the asymptotes close in infinity, suggesting a rich topological and symmetrical structure. This directly connects with topological compacting concepts, global asymptotic symmetries, and stable field solutions, suggesting possible connections between geometric structures and asymptotic field-like behaviors.

\section{Conclusion}

The use of the system of ordinary differential equations allowed important characteristics to be mapped, revealing a dynamic with the existence of nodes. In these nodes, the topology changes when the rotation between the planes in the division occurs. Zero rotation is associated with Euclidean space (locally). The simulation showed convergence and was dynamically stable close to $0 + n\pi$ radians. Theoretically, Euclidean space (locally) has an ODE system that represents the trajectory of this stable and convergent dynamic. However, the rotation close to $\pi/2 + n\pi$ radians, associated with the Alpha group, suggests asymptotic attractor-like behavior, indicating the existence of an attractor that converges asymptotically at infinity. Matrix $A$  corresponds to a regime of maximal topological deformation; its interpretation suggests that as $\theta$ approaches $\pi/2$, the nature of the transformations represented becomes more complex and may involve a significant deformation in the space topology (global). Briefly, the transformation properties induced by the Alpha Group and its matrix $A$ exhibit structural features analogous to gauge-like symmetries acting on geometric configurations. The matrix $A$ may therefore be interpreted as generating structured symmetry-like transformations in $R^4$. The geometric point ($\mu$) defined in the alpha group is a way of unifying concepts of singularity with stable geometric structures.
A key dynamical and topological outcome is the transition of topological domains from a collapsed 3-sphere ($S^3$), which corresponds to a real regime, to a 4-sphere-like regime ($S^4$), representing a complex and hypercomplex regime. In this transition, a hidden degree of freedom begins to emerge through the imaginary component of the structure.

This process is observed in the dynamical behavior of the complex eigenvalues of $M(\theta)$, which play a central and significant role. When the real part $\operatorname{Re}(\lambda) \rightarrow 0$ and the imaginary part $\operatorname{Im}(\lambda) \rightarrow \infty$, phenomena such as may indicate torsion-like or curvature-like effects, may be associated with resonance phenomena and Hopf-like bifurcation behavior.

In this context, the vector $\mu$ acts as an idempotent and invariant element, ensuring internal coherence and stability throughout the topological deformations. Its presence is essential for preserving the structural integrity and algebraic consistency of the Alpha Group framework.

Under these conditions, $\mu$ generates a type of attractor that, while dynamically stable, asymptotically converges at infinity. This results in a form of internal dynamic topological coherence, characteristic of the transition between real and complex topological domains.

\end{document}